\documentclass[12pt]{amsart}
\usepackage{amssymb}
\usepackage{xspace}

\newtheorem{Thm}{Theorem}[section]

\newtheorem{Lem}[Thm]{Lemma}
\newtheorem{Prop}[Thm]{Proposition}

\def\ldots{\mathinner{\ldotp\ldotp\ldotp}}
\def\ldots{\mathinner{\cdotp\cdotp\cdotp}}

\def \cal{\mathcal}

\newcommand\ds{\displaystyle}
\newcommand\LB[1]{\label{#1}}
\newcommand\BE[2]{\begin{#1} #2 \end{#1}}
\newcommand\EQ[2]{\BE{equation}{\LB{#1} #2}}

\newcommand\ZZ{{\mbox{\sf Z\kern-.45em Z}}}
\newcommand\vv{\kern.344em{\rule[.18ex]{.075em}{1.32ex}}\kern-.344em}

\newcommand\RE{{\mbox{\rm I\kern-.21em R}}}
\newcommand\NN{\mbox{I\kern-.21em N}}

\newcommand\Rb{Riesz basis\xspace}

\newcommand\op{operator\xspace}

\newcommand\Hs{\ensuremath{H^s(-\pi,\pi)}\xspace}

\newcommand\E{\ensuremath{\mathcal E}\xspace}

\newcommand\Sss{ Sobolev spaces\xspace}

\newcommand\iy{\infty}

\newcommand\PI{(-\pi,\pi)}

\def\th{\theta}

\def\a{\alpha}

\newcommand\en{e_{\lambda_n}}

\newcommand\dl{\delta}

\newcommand\lmn{\lambda_n}
\newcommand\Lm{\Lambda}

\def\l{\left} \def\r{\right}

\begin{document}

\title
{Interpolation of subspaces and applications to exponential bases
in Sobolev spaces}
 \pagestyle{plain}
\author{Sergei Ivanov}
\address{Russian Center of Laser Physics,
St.Petersburg State University, Ul'yanovskaya 1,198904 St.Petersburg, Russia}
\email {Sergei.Ivanov@pobox.spbu.ru}
\author{Nigel Kalton}
\address{Department of Mathematics
University of Missouri
Columbia
Mo. 65211 USA}
\email {nigel@math.missouri.edu}
\thanks{
The first author  was partially supported by the Russian Basic
Research Foundation
(grant \# 99-01-00744) and the second author was supported by NSF-grant
DMS-9870027}
\thanks{
The first author is grateful to
S.Avdonin for  fruitful discussions.}

 \begin{abstract}
We give precise conditions under which the real interpolation space
$[Y_0,X_1]_{\theta,p}$ coincides with a closed subspace of
$[X_0,X_1]_{\theta,p}$ when $Y_0$ is a closed subspace of codimension one.
We then apply this result to nonharmonic Fourier series in Sobolev spaces
$H^s\PI$ when $0<s<1.$ The main result: let \E be a family of exponentials
$\exp(i \lambda_n t)$ and \E forms an unconditional basis in $L^2\PI$.
Then there exist two number $s_0, s_1$ such that \E forms an
unconditional basis in $H^s$ for $s<s_0$, \E forms an unconditional basis
in its span with codimension 1 in $H^s$ for $s_1<s$. For $s_0\le s \le s_1$
the exponential family is not an unconditional basis in its span.
\end{abstract}

\maketitle

\section{Introduction}\LB{n0}

In this paper we will apply a result on interpolation of
subspaces
 to the study of exponential Riesz bases in Sobolev spaces.

In  section \ref{interpolation} we consider the comparison of  the
 interpolation  spaces
$X_\th:=[X_0,X_1]_{\th,p}$ and
$Y_\th:=[Y_0,X_1]_{\th,p}$ for $1\le p<\infty$, where $Y_0$  is a
subspace of
$X_0$
with codimension one, say $Y_0=\text{ker }\psi$ where $\psi\in X_0^*.$
 This problem, as far as we know, was first formulated
in
\cite{LM},
v.1, Ch.1.18 in 1968. As we show in Theorem \ref{interpol1} there are two
indices $0\le \sigma_0\le \sigma_1\le 1$ which may be explicitly
evaluated in terms of the $K$-functional of $\psi$ so that:
\begin{enumerate} \item If
$0<\theta<\sigma_0$ then
$Y_{\theta}$ is a closed subspace of codimension one in $X_{\theta}$.
\item If
$\sigma_1<\theta<1$ then $Y_{\theta}=X_{\theta}$ with equivalence of
norm
and \item If $\sigma_0\le \theta\le\sigma_1$ then the norm on
$Y_{\theta}$ is not equivalent to the norm on $X_{\theta}.$
\end{enumerate}

Let us discuss the history of this theorem.
 The special case of a Hilbert space of Sobolev
type connected with elliptical boundary data was considered in \cite{LM},
and in this case the critical indices $\sigma_0$ and $\sigma_1$
coincide.
In the well
known case \cite{LM}
$X_1=L^2(0,\iy)$,
$X_0=W_2^1(0, \iy)$ and
$Y_0$ is the subspace of $W_2^1$ of functions vanishing at the
origin, this critical value is
$\sigma_0=\sigma_1=1/2.$  Later R. Wallsten \cite{Wallsten} gave an
example where the critical indices satisfy $\sigma_0<\sigma_1.$
The general problem  was considered by J. L\"ofstr\"om
\cite{L1}, where some special cases of Theorem \ref{interpol1} are
obtained. Later, L\"ofstr\"om in an unpublished (but web-posted) preprint
from 1997, obtained most of the conclusion of Theorem \ref{interpol1}:
specifically he obtained the same result except he did not treat the
critical values $\theta=\sigma_0,\sigma_1.$  The authors were not aware
of L\"ofstr\"om's earlier work during the initial preparation of this
article and our approach is rather different.    A  more general
but closely related problem on interpolating subspaces of codimension
one has been recently considered in \cite{KMP} and \cite{KS}.
For general results on subcouples we refer to
 \cite{Janson}.

Let us recall next that a sequence $(e_n)_{n \in\mathbb Z}$ in a Hilbert
space $\mathcal H $ is called a {\it Riesz basic sequence} if there is a
constant $C$ so that for any finitely non-zero sequence
$(a_n)_{n\in\mathbb Z}$ we have
$$ \frac1C \left(\sum_{n\in\mathbb Z}|a_n|^2\right)^{\frac12} \le
\|\sum_{n\in\mathbb Z}a_ne_n\| \le C
 \left(\sum_{n\in\mathbb Z}|a_n|^2\right)^{\frac12}.$$
A {\it Riesz basis} for $\mathcal H$ is a Riesz basic sequence whose
closed linear span $[e_n]_{n\in\mathbb Z}=\mathcal H.$  A sequence
$(e_n)$ is an {\it unconditional basis,}  respectively {\it unconditional
basic sequence} if
$(e_n/\|e_n\|)_{n\in\mathbb Z}$ is a Riesz basis, respectively, a Riesz
basic sequence.

In the second part of the paper we apply our interpolation result
to study the basis properties of exponential families
$\lbrace {e^{i\lambda_n t} \rbrace} $ in \Sss.
These families appear
in such fields of mathematics as the theory of dissipative
operators (the Sz.--Nagy--Foias model)  ,
the Regge problem for
resonance scattering, the theory of initial boundary value
problems,
control theory for distributed parameter systems, and
signal processing, see, e.g.,
\cite{Nik}, \cite{IP78},
\cite{LLT}, \cite{AI95}, \cite{Seip95}.
One of the most important problems arising in all
of these applications is the question of the Riesz basis
property od these families.
In the space $L^2\PI$ this problem
has been studied
for the first time in the classical work of Paley
and Wiener  \cite{PW34}.
The problem has now a complete solution
\cite {KNP},  \cite {Minkin} on the basis of an approach suggested by
B.~S.~Pavlov.

The principal result for Riesz bases  can be formulated as
follows
\cite{KNP}.
\begin{Prop}\LB{A2}
 The sequence $(e^{i\lambda_n t})_{n\in\mathbb Z}$ is a Riesz
basis for $L^2\PI$ if and only if $\sup|\Im \lambda_n|<\infty,$
\EQ{2.1}{
\inf_{k\ne j}\left|\lambda_k-\lambda_j\right|>0.
} and there is an entire function $F$ of exponential type $\pi$ (the {\it
generating function}) with simple zeros at $(\lambda_n)_{n\in\mathbb
Z} $   and  such that for some $y$
$\left|F(x+iy)\right|^2$ satisfies the Muckenhoupt condition $(A_2)$
(we shall write  this as $|F|^2\in(A_2)$):
$$
\sup_{I\in {\cal J}}\Bigl\{
\frac{1}{\left|I\right|}
\int_I\left|F(x+iy)\right|^2\,dx\,\,\frac{1}{\left|I\right|}
\int_I\left|F(x+iy)\right|^{-2}\,dx\Bigr\}\,\, < \,\,\infty,
$$
where ${\cal J}$ is the set of all intervals of the real axis.
\end{Prop}

In \cite{Minkin} a corresponding characterization is given for exponential
families which form an unconditional basis of $L^2(-\pi,\pi)$ when
$\Im\lambda_n$ can be unbounded both from above and below.

\par
Let us describe known results concerning
exponential bases in Sobolev spaces.
The first result in this direction has been obtained
by D.~L.~Russell in \cite{Russell82}
Russell studied the unconditional basis property for exponential
families in the \Sss $H^m (-\pi,\pi)$ with  $m\in \mathbb Z$.
\par
\begin{Prop}\LB{Russell}
\cite{Russell82} Suppose $(e^{i\lambda_nt})_{n\in\mathbb Z}$ is a Riesz
basis for $L^2\PI.$ Suppose $m\in\mathbb N$ and suppose $\mu_1,\ldots,
\mu_m\in\mathbb C\setminus \{\lambda_n: \ n\in\mathbb Z\}$ are distinct.
Then $(e^{i\lambda_nt})_{n\in\mathbb Z}\cup (e^{i\mu_kt})_{k=1}^m$ is an
unconditional basis of $H^m\PI.$   In particular
$(e^{i\lambda_nt})_{n\in\mathbb Z}$ is an unconditional basic sequence
whose closed linear span has codimension $m$ in $H^m\PI.$\end{Prop}
\par

\par In \cite{NS86} the unconditional basis property for an exponential
family was studied
in  $H^s(-\pi,\pi)$  for noninteger $s$ for the case
$\lambda_n$ being the eigenvalues of a Sturm--Liouville operator with
a smooth potential.

Note that the
generalization of the  Levin--Golovin theorem for Sobolev
spaces has been obtained \cite{AI00}
using `classical methods' of the entire function theory. Suppose
$\{ \lambda_n \}_{n\in \mathbb Z}$ are the zeros
of an entire function $F$ of exponential type $\pi$,
$(\lambda_n)$ is separated (\ref{2.1}), and
on some line  $\{ x+iy \}_{x\in \RE}$ we have
$$  C^{-1}(1+|x|)^s\le
\vert F(x+iy) \vert \le C(1+\vert x \vert)^s.
$$
Then
the family $\{e^{i\lambda_kt}/(1+|\lambda_k|)^s\}$ forms a Riesz basis
in $H^s\PI$.
Notice that this result was applied to several
controllability problems for the wave type equation \cite{AIR}.

Recently Yu.~Lyubarskii and K.~Seip \cite{LS} have established a
necessary and sufficient
  criterion for
sampling/interpolation problem for weighted Paley-Wiener spaces, which
gives a criterion for a sequence to be an unconditional basis in $H^s.$
 For the case, $\sup |\Im
\lambda_n|<\infty,$ the main result is the following:

\begin{Thm}\label{LS}
 $(e^{i\lambda_nt})_{n\in\mathbb Z}$  forms an unconditional basis
 in
\Hs
if and only if $(\lambda_n)$ is
separated (i.e. (\ref{2.1}) holds) and for the generating function $F$ we
have
$|F(x+iy)|^2/(1+|x|^{2s})\in (A_2)$ for some $y$.\end{Thm}

The main idea of the present paper is that if
$(e^{i\lambda_nt})_{n\in\mathbb
Z}$ forms a Riesz basis in $L^2(-\pi,\pi)$ then it also forms an
unconditional basis of a subspace $Y_0$ of $H^1(-\pi,\pi)$ of codimension
one.  Then, by interpolation, one obtains that
$(e^{i\lambda_nt})_{n\in\mathbb Z}$ is an unconditional basis of the
intermediate spaces $[Y_0,L^2]_{\theta,2}$ for $0<\theta<1.$ This
approach was suggested in
\cite{I96} by the first author.
The main result of \cite{I96} is incorrect in the general case
because a mistake
connected with interpolation of subspaces. Here we correct this mistake.

Let us describe the results concerning unconditional bases in \Sss.
One of our
main results for Riesz bases is as follows:

\begin{Thm}\label{main}  Suppose $(e^{i\lambda_nt})_{n\in\mathbb Z}$
forms a Riesz basis of $L^2(-\pi,\pi).$
Suppose $(\lambda_n-n)_{n\in\mathbb Z}$ is bounded and
let $\delta_n=\Re\lambda_n-n.$
Then there exist critical indices
$0<s_0\le s_1<1$ given by:$$ s_1= \frac12-\lim_{\tau\to\infty}
\inf_{t\ge
1}\frac{1}{\log\tau}\sum_{t< |n|\le\tau
t}\frac{\delta_n}{n}$$ and
$$s_0=\frac12 - \lim_{\tau\to\infty}
\sup_{t\ge
1}\frac{1}{\log\tau}\sum_{t<|n|\le\tau
t}\frac{\delta_n}{n}$$ such that: \newline
(1)
$(e_{\lambda_n})_{n\in\mathbb Z}$ is an unconditional  basis
of  the  Sobolev  space   $H^s$  if  and  only   if
$0\le s< s_0.$
\newline
(2)  $(e_{\lambda_n})_{n\in\mathbb  Z}$  is an
unconditional  basis  of  a  closed  subspace  of  $H^s$  of
codimension one if and only if $s_1<s\le 1.$\newline (3)  If
$s_0\le   s\le   s_1$   then   $(e_{\lambda_n})$  is  not  an
unconditional basic sequence. \end{Thm}

This result is deduced from results in Sections \ref{n2} and \ref{n3}.
In Section \ref{n3} we in fact consider the more general situation for
unconditional bases and give rather more technical results.  The above
Theorem \ref{main} however is the simplest case and follows by combining
Theorem \ref{uncbasis}, Theorem \ref{perturbation3} and  Theorem
\ref{strictineq}.  Our approach is based on estimates of the
K-functional for the continuous linear functional on $H^1(-\pi,\pi)$
which annihilates each $e^{i\lambda_nx}$ whose existence is guaranteed by
the result of Russell (Proposition \ref{Russell}).  The estimates are in
terms of the generating function $F.$

Once one has Theorem \ref{main} then it is easy to construct real
sequences $(\lambda_n)$ to show that $s_0,s_1$ can take any values in
$(0,1)$ such that $s_0\le s_1.$  In the case of regular power behavior of
$F$
i.e. for some $y\ge 0,$ $|F(x+iy)|\sim (1+|x|)^{s}$  one has
$s_1=s_0=s+\frac12.$

The results for the whole scale \Hs can then be obtained by
`shift' using the fact that  the differentiation \op with
appropriate conditions is an isomorphism between a
one-codimensional subspace of
$H^m$ and $H^{m-1}$; we will not pursue this extension.

\section{Interpolation of
subspaces}\label{interpolation}\setcounter{equation}{0}

Let $(X_0,X_1)$ be a Banach couple with $X_0\cap X_1$  dense
in $X_0,X_1.$ If $0<\theta<1$  and $1\le p<\infty$ the  real
interpolation  space  $X_{\theta}=[X_0,X_1]_{\theta,p}$   is
defined, see, e.g., \cite{BS},
to  be  the  set  of  all  $x\in X_0+X_1$ such that
\begin{equation*}\label{interdef}     \|x\|_{X_\theta}=\left(
\int_0^{\infty}t^{\theta                        p-1}K(t,x)^p
dt\right)^{\frac1p}<\infty,\end{equation*}
where $K(t,x)$ is the $K$--functional.
An    equivalent
definition  \cite{BS} p. 314  (yielding  an  equivalent  norm) can be
given by using   the J-method:
\begin{equation*}
\label{interdef2}
\|x\|_{X_{\theta}}=\inf\left\{\left(\sum_{k\in\mathbb
Z}\max\{\|x_k\|_0,2^k\|x_k\|_1\}^p\right)^{\frac1p}:\
x=\sum_{k\in\mathbb Z}2^{\theta k}x_k\right\},
\end{equation*}
where the series converges in $X_0+X_1.$

Now  suppose  $0\neq  \psi\in  X_0^*$  and  let $Y_0$ be its
kernel. We suppose also (only this case is interesting)
that $Y_0\cap X_1$ is dense in $X_1$, i.e., $\psi$ is not bounded
in $X_1$.

Let  $Y_{\theta}$  be  the  corresponding  spaces
obtained   by   interpolating   $Y_0$   and  $X_1.$  Clearly
$Y_{\theta}\subset X_{\theta}$  and the  inclusion has  norm
one.  It is easy to show that the closure of $Y_{\theta}$ in
$X_{\theta}$ is  either a  subspace of  codimension one when
$\psi$  is  continuous  on  $X_{\theta}$  or  the  whole  of
$X_{\theta}$ when $\psi$ is not continuous.

Let    us    now    introduce    two    important   indices.
\begin{equation*}\LB{s001}                           \sigma_1
=\lim_{\tau\to\infty}\sup_{0<\tau  t\le  1} \frac1{\log\tau}
\log\frac{K(\tau t,\psi)}{K(t,\psi)}
\end{equation*}
and
\begin{equation*}\label{s002}
\sigma_0=\lim_{\tau\to\infty}\inf_{0<\tau  t\le
1}\frac{1}{\log\tau}
\log\frac{K(\tau t,\psi)}{K(t,\psi)},
\end{equation*}
where
$
K(t,\psi)=K(t,\psi; X_0^*,X_1^*).
$
>From the multiplicative properties of the function
$K(\tau t,\psi)/K(t,\psi)$
it  is clear  that  these   limits  exist  and   $0\le  \sigma_0\le
\sigma_1\le  1.$  Since  $K(t,\psi)$  is  bounded   as
$t\to\infty$   we   can    also   write:
$$   \sigma_1
=\lim_{\tau\to\infty}\sup_{0<t<\infty}      \frac1{\log\tau}
\log\frac{K(\tau t,\psi)}{K(t,\psi)}.
$$

Let  us  observe  that:
\begin{equation*}\label{Kdef} \sup
\{|\psi(x)|:         \    \max\{\|x\|_0,t\|x\|_1\}\le     1\}=
K(t^{-1},\psi)  \end{equation*}
We  define  a  sequence
$(w_n)_{n\in\mathbb Z}$ by  $$ w_n =  K(2^{-n},\psi)^{-1}.$$
Notice        that        $\inf_{n\in\mathbb        Z}w_n\ge
\|\psi\|_{X_0^*}^{-1}>0$   and   that   in   general  $w_n\le
w_{n+1}\le   2w_n.$   Now   it   is   easy   to   see   that
\begin{eqnarray*}\LB{s01} \sigma_1  &=\lim_{k\to\infty}\sup_n
\frac1k                          \log_2\frac{w_{n+k}}{w_n}\\
\sigma_0&=\lim_{k\to\infty}\inf_{n\ge          0}\frac{1}{k}
\log_2\frac{w_{n+k}}{w_n}.   \end{eqnarray*}

As mentioned in the introduction the following result is a slight
improvement of a result of L\"ofstr\"om \cite{L2}, who obtains the same
result by quite different arguments  except for the critical indices
$\theta=\sigma_0,\sigma_1.$

\begin{Thm}\LB{interpol1}  1.  $Y_{\theta}=X_{\theta}$ (with
equivalence    of   norm)     if    and    only     if    $
\theta>\sigma_1$.\newline
2.   $Y_{\theta}$  is   a  closed
subspace of codimension one  in $X_{\theta}$ if and  only if
$\theta<\sigma_0.$\newline
3.If   $\sigma_0\le   \theta\le
\sigma_1$ then $Y_{\theta}$  is not closed  in $X_{\theta}.$
\end{Thm}

We shall consider the weighted $\ell_p$ space $\ell_p(w)$ of
all  sequences  $(\alpha_n)_{n\in\mathbb  Z}$  such  that $$
\|\alpha\|=                          \left(\sum_{k\in\mathbb
Z}w_n^p|\alpha_n|^p\right)^{\frac1p}.$$
We  shall use  $\zeta_n$
for the standard basis vectors.  On $\ell_p(w)$ we  consider
the shift operator $S((\alpha_n))= (\alpha_{n-1}).$ From the
above remarks it is  clear that $S,S^{-1}$ are  both bounded
and  $\|S\|\le  2,\|S^{-1}\|=  1.$  Furthermore the spectral
radius  formula  shows  that  $2^{\sigma_1}$ is the spectral
radius  $r(S)$  of  $S.$  Now  let  $P_+$  be the projection
$P_+(\alpha)=  (\delta_n\alpha_n)$  where  $\delta_n=1$   if
$n\ge  0$  and  $0$  otherwise.    It  is  easy to calculate
\begin{equation*}
\label{norm}   \|P_+S^{-n}\|=   \sup_{k\ge
0}\frac{w_k}{w_{n+k}}
\end{equation*}
and so this implies that
$r(P_+S^{-1})=2^{-\sigma_0}.$

We     will     need     the     following     key    Lemma:
\begin{Lem}\label{interpol2}   Let   $0<\theta<1$   and  let
$T_{\theta}=S-2^{\theta}I.$ Then \begin{enumerate}
\item $T_\theta$
is  an   isomorphism  onto   $\ell_p(w)$  if   and  only  if
$\sigma_1<\theta.$
\item $T$ is an isomorphism onto a proper
closed subspace  if and  only if  $\theta<\sigma_0.$ In this
case the range of $T$ is the subspace of codimension one  of
all      $\alpha$      such      that     $\sum_{n\in\mathbb
Z}2^{n\theta}\alpha_n=0.$\end{enumerate}
\end{Lem}

\begin{proof} First observe  that if $\theta>\sigma_1$  then
$T_{\theta}$ must be  an isomorphism onto  $\ell_p(w)$ since
$2^{\theta}$ exceeds the spectral radius of $S.$ Furthermore
since the spectrum of $S$ is invariant under rotations it is
clear  that  $T_{\sigma_1}$  cannot  be  an isomorphism onto
$\ell_p(w).$ Also note that $T_{\theta}$ is always injective
and that  if $f_\th$  is a  linear functional annihilating
its range then $f_\th(\zeta_n)=c2^{n\theta}$ for some constant
$c$, i.e.,
$f_{\theta}(\alpha)=\sum_{n\in\mathbb Z}2^{n\theta}\alpha_n=0$.
 This implies that  the closure of the range  is either
the  whole  space  or  the  subspace  of codimension one
when  $\sum_{n\in\mathbb
Z}2^{n\theta q}w_n^{-q}<\infty$.   Here  $\frac1p+\frac1q=1$
and the formula must be modified if $p=1.$

We next show that if $\theta<\sigma_0$ then $T_{\theta}$  is
an isomorphism onto a closed subspace of codimension one.

Next let  $E=[\{\zeta_n:   \ n\le  -1\}]$ and  $F= [\{\zeta_n:\ n\ge
1\}].$ We remark that  $T_{\theta}(E)$ is easily seen  to be
closed because  $T_{\theta}$ is  an isomorphism  on the {\it
unweighted }$\ell_p$ and $w_n$ is bounded for $n\le -1.$  If
we show $T_{\theta}(F)$ is closed then we are done, since it
is clear  this will  imply that  $T_{\theta}(E+F)$ is closed
and this  is a  subspace of  co-dimension one  in the range.
However       $2^{-\theta}>r(P_+S^{-1})$       so       that
$2^{-\theta}-P_+S^{-1}$ is an  isomorphism.  Restricting  to
$F$   this   implies   $(2^{-\theta}-S^{-1})F$   and   hence
$T_{\theta}(F)$ is closed.

The  proof  is  completed  by  showing  that  if  $\theta\le
\sigma_1$ then if $T_{\theta}$  if has closed range  it must
satisfy $\theta<\sigma_0.$ Note first  that it is enough  to
establish  this  for  $\theta<\sigma_1$  since  the  set  of
operators  with  Fredholm  index  one  is  open.     Suppose
$\sigma_0<\theta<\sigma_1$ and $T_{\theta}$ is closed.  Then
$T_{\theta}$ has  a lower  estimate $\|T_{\theta}\alpha\|\ge
c\|\alpha\|$   for   all   $\alpha$   where   $c>0.$  Assume
$w_{n+k}>2^{n\theta}w_k$  for   some  $n\in\mathbb   N$  and
$k\in\mathbb      Z.$      Then      consider       $\alpha=
(I+2^{-\theta}S+\cdots+2^{-n\theta}S^{n})^2\zeta_k.$  Note  that
$\|\alpha\|\ge        n2^{-n\theta}w_{n+k}.$         However
\begin{align*}
\|T_{\theta}^2\alpha\|=&
2^{2\th}w_k + 2\cdot2^{(-n+1)\theta}w_{n+k+1}+2^{-2n\theta}w_{2n+k+2}
\\ &\le
8\max\{w_k,2^{-n\theta}w_{n+k},2^{-2n\theta}w_{2n+k}\}.    \end{align*}
Let  $v_n=  2^{-n\theta}w_n.$  Then  we  have  if $nc^2>8,$ $$
(nc^{2}-8)v_{k+n}   \le  8\max\{v_k,v_{k+2n}\}.$$    In
particular if $nc^2>16,$ $$ v_{k+n}< \max\{v_k,v_{k+2n}\}.$$

Now  since  $\theta<\sigma_1,$  we  can find
$k\in\mathbb Z,n>16c^{-2}$
so    that    $w_{n+k}<2^{n\theta}w_k$     or
$v_{n+k}<v_k.$      Iterating       gives      us       that
$(v_{k+rn})_{r=0}^{\infty}$ is monotone increasing.
Now for
any   large   $N$   and   any   $j\ge   0$   we   have
$$
\frac{w_{j+N}}{w_j}   \ge   \frac{w_{k+r_2n}}{w_{k+r_1n}}\ge
2^{n(r_2-r_1)\theta}$$
where $r_1$, $r_2$  such that
$k+(r_1-1)n\le j\le k+r_1n $ and
$k+r_2n\le j+N\le k+(r_2+1)n$.
This gives us
$$   \frac{w_{j+N}}{w_j}\ge 2^{(N-2n)\theta}.$$
Hence  $$  \inf_{j\ge  0}   \frac1N
\log_2 \frac{w_{j+N}}{w_j} \ge (1-\frac{2n}{N})\theta.$$
Letting  $N\to\infty$  gives  $\sigma_0\ge  \theta.$ To show
that in  fact $\theta<\sigma_0$  needs only  the observation
again  that  the  set  of  $\theta$  where  $T_{\theta}$ has
Fredholm index one is open.\end{proof}

We  now  use  Lemma  \ref{interpol2}  to  establish our main
result Theorem \ref{interpol1} on interpolating subspaces:

\begin{proof}   Let   us   suppose   next  that  either  (a)
$\theta<\sigma_0$  or  (b)  $\theta>\sigma_1.$  This implies
there  exists   a  constant   $D$  so   that  $\|\alpha\|\le
D\|T_{\theta}\alpha\|$ for all $\alpha\in\ell_p(w)$; in case
(a)  $T_{\theta}$  maps  onto  the  subspace of $\ell_p(w)$
defined       by       $f_{\theta}(\alpha)=\sum_{n\in\mathbb
Z}2^{n\theta}\alpha_n=0$, while in case (b) $T_{\theta}$  is
an   isomorphism   onto   the   whole   space   (see   Lemma
\ref{interpol2}).      We   observe   that   in   case    (a)
 the linear functional $\psi$ extends  to
a  continuous  linear  functional  on  $X_{\theta}$  as   $$
\sum_{n\in\mathbb Z} 2^{n\theta}K(2^n,\psi)<\infty.$$

Now  suppose  $x\in  X_{\theta}$  with $\|x\|_{X_{\theta}}=1$
with the additional assumption in case (a) that $\psi(x)=0.$
Then  we   may  find   $(x_n)_{n\in\mathbb  Z}$   such  that
$\sum_{n\in\mathbb Z}2^{\theta n}x_n=x$  and $$  \left(\sum_{k\in\mathbb
Z}\max\{\|x_k\|_0,2^k\|x_k\|_1\}^p\right)^{\frac1p}\le    2.$$
Then
$$
\left(\sum_{n\in\mathbb Z}|\psi(x_n)|^pw_n^p\right)^{\frac1p}\le 2,
$$
since
$$
|\psi(x)|\le w_n^{-1} \max\{\|x\|_0,2^n\|x\|_1\}.
$$
In case (a) we
additionally            have             $$\sum_{n\in\mathbb
Z}2^{n\theta}\psi(x_n)=0.$$     Thus     we     can     find
$\alpha\in\ell_p(w)$ with $  T_{\theta}(\alpha)=(\psi(x_n))$
and $\|\alpha\|\le  2D.$ Then  we can  find $u_n\in  X_0\cap
X_1$    such    that
$\max\{\|u_n\|_0, 2^n\|u_n\|_1\}\le2|\alpha_n|w_n$
and   $\psi(u_n)=\alpha_n.$   Let    $v_n=
u_{n-1}-2^{\theta}u_{n}.$  Then  $$  \left(\sum_{k\in\mathbb
Z}\max\{\|v_k\|,2^k\|v_k\|_1\}^p\right)^{\frac1p}\le       16D
\|x\|_{X_{\theta}}.$$                                    Now
$\psi(v_n)=\alpha_{n-1}-2^{\theta}\alpha_n=\psi(x_n)$    and
$\sum_{n\in\mathbb   Z}2^{n\theta}v_n=0.$   Hence
$$
x=\sum_{n\in\mathbb  Z}2^{\th n}(x_n-v_n)
$$
and  so  $x\in Y_{\theta}$
with  $\|x\|_{Y_\theta}\le  (16D+2)\|x\|_{X_{\theta}}.$ From
this it follows that  in case (a) we  have $Y_{\theta}=\{x:\
\psi(x)=0,\   x\in   X_{\theta}\}$    and   in   case    (b)
$Y_{\theta}=X_{\theta}.$

Next we  consider the  converse directions.   Assume  either
(aa)   $\psi$    is   continuous    on   $X_{\theta}$    and
$Y_{\theta}=\{x:\  \psi(x)=0,\  x\in  X_{\theta}\}$  or (bb)
$Y_{\theta}=X_{\theta}.$ In either case there is a  constant
$D$ so that if $x\in Y_{\theta}$ then $\|x\|_{Y_{\theta}}\le
D\|x\|_{X_{\theta}}.$ Observe  that in  case (a)  the linear
functional $f_{\theta}$ is continuous on $\ell_p(w)$ and  so
the range  of $T_{\theta}$  is contained  in its  kernel; in
case (bb) its range is dense.

Assume $\alpha=(\alpha_n)_{n\in\mathbb Z}\in \ell_p(w)$ with
$\|\alpha\|=1;$    in    case    (aa)    we    also   assume
$f_{\theta}(\alpha)=0.$ We first find $x_n\in X_0\cap  X_1$
with    $\psi(x_n)=\alpha_n$    and    so    that
$\max\{\|x_n\|_0,2^n\|x_n\|\}\le  2|\alpha_n |w_n$
for $n\in\mathbb
Z.$ Let $x=\sum_{n\in\mathbb Z}2^{n\theta}x_n$ so that $x\in
X_{\theta}$ with $\|x\|_{X_{\theta}}\le  2.$ In case  (aa) we
have additionally  that $\psi(x)=f_{\theta}(\alpha)=0.$  Now
we can find $y_n\in Y_0\cap X_1$ so that  $\sum_{n\in\mathbb
Z}2^{n\theta}y_n=x$    and
$$    \left(\sum_{k\in\mathbb
Z}\max\{\|y_k\|,2^k\|y_k\|_1\}^p\right)^{\frac1p}\le 4D.$$

Now             let             $u_n=x_n-y_n$            and
$v_n=\sum_{k=n+1}^{\infty}2^{(k-n-1)\theta}u_k.$   Then
$$\left(\sum_{k\in\mathbb
Z}\max\{\|u_k\|,2^k\|u_k\|_1\}^p\right\}^{\frac1p}\le   4D+2.$$
We argue that
\begin{equation}\label{useful}
\left(\sum_{k\in\mathbb
Z}\max\{\|v_k\|_0,2^k\|v_k\|_1\}^p\right)^{\frac1p}\le
C_{\theta}(4D+2),\end{equation}                        where
$$
C_{\theta}=(\sum_{k<0}2^{k\theta}+\sum_{k\ge
0}2^{k(\theta-1)}).$$
To  show (\ref{useful})  we note  that
$$
2^n\|v_n\|_1                    \le
\sum_{k=n+1}^{\infty}2^{(k-n-1)(\theta-1)}2^k\|u_k\|_1
$$
and (since $\sum 2^{n\th}u_n=0$)
$$                       \|v_n\|_0                       \le
\sum_{k=-\infty}^{n}2^{(k-n-1)\theta}\|u_k\|_0.$$

Let   $\beta_n=\psi(v_n).$   Then   $\beta\in\ell_p(w)$  and
$\|\beta\|\le        C_{\theta}(4D+2).$        But       now
$(T_{\theta}(\beta))_n   =\psi(u_n)=\psi(x_n)=\alpha_n$   so
that  $T_{\theta}$  is  an  isomorphism  onto  the kernel of
$f_{\theta}$ in case (aa)  or onto $\ell_p(w)$ in case (bb).
These two cases  combined with the observation  that $Y_{\theta}$
can  only  be  a  proper  closed subspace of $X_{\theta}$ if
$\psi$ is continuous on  $X_{\theta}$ complete the proof  of
the Theorem.
\end{proof}

\section{Sobolev spaces}\LB{n2}
\setcounter{equation}{0}

In this section we investigate a special case of the results of the
previous section for Sobolev spaces. These results are preparatory for
Section \ref{n3} where we apply them to exponential bases. Let
$L^2=L^2(-\pi,\pi)$
and let us denote the standard inner-product on $L^2(-\pi,\pi)$ by
$$ (f,g)=  \int_{-\pi}^{\pi}f(x)\overline{g(x)}dx.$$  We denote by
$\|f\|$ the standard norm on $L^2.$

 For $s>0$  we define the Sobolev
space
$H^s(\mathbb R)$ to be the space of  all $f\in L_2(\mathbb R)$   so
that
$$   \|f\|^2_{H^s}   :=  \int_{-\infty}^{\infty}|\hat
f(\xi)|^2(1+|\xi|^{2s})d\xi<\infty$$
($\hat f$ is the Fourier transform).
We   then   define  the
Sobolev  space  $H^s=H^s(-\pi,\pi)$  to  be  the  space   of
restrictions of $H^s(\mathbb  R)-$functions to the  interval
$(-\pi,\pi)$ (with the obvious induced quotient norm).  When
$s=1$  the  space  $H^1$  reduces  to  the  space  of $f\in
L^2(-\pi,\pi)$ so  that $f'\in  L^2$ under  the (equivalent)
norm:            $$      \|f\|_1^2     =     \int_{-\pi}^{\pi}
|f(t)|^2+|f'(t)|^2dt<\infty.$$ Then if $0<s<1$ we have $H^s=
[H^1,L^2]_{1-s}=[H^1,L^2]_{1-s,2}$ \cite{LM}.

For   $z\in   \mathbb   C$   we   define  $e_z(x)=e^{izx}\in
L^2(-\pi,\pi).$  Now suppose $\psi\in (H^1)^*$; we define its Fourier
transform $F=\hat\psi$ to be the entire function $F(z):=\psi(e_z)$ for
$z\in
\mathbb C.$ Let us first identify $(H^1)^*$ via its Fourier transform:

\begin{Prop}\label{dualsobolev} Let $F$ be an entire function.  In order
that there exists $\psi\in (H^1)^*$ with $F=\hat\psi$ it is necessary and
sufficient that:
\begin{equation}\label{type} F \text{ is of exponential type }\le \pi.
\end{equation}
\begin{equation}\label{integral}
\int_{-\infty}^{\infty}\frac{|F(x)|^2}{1+x^2}dx<\infty.
\end{equation} These conditions imply the estimate:
\begin{equation}\label{type2} \sup_{z\in\mathbb C}
\frac{|F(z)|}{(1+|z|)e^{\pi |\Im z|}}<\infty.\end{equation} \end{Prop}

\begin{proof} These results follow immediately from the Paley-Wiener
theorem once one observes that $\psi\in (H^1)^*$ if and only if $\psi$ is
of the form $$\psi(f) = \alpha f(0) + \varphi(f')$$ where $\varphi\in
(L^2)^*.$
\end{proof}

Consider  $H^1$  with  the   inner  product:    $$   \langle
f,g\rangle_t  =  (f',g')+t^2  (f,g)$$ where $t>0.$
Let us denote by $\|\psi\|_t$ the norm of $\psi$ with
respect to $\|\cdot\|_t$ where
$\|f\|_t^2=\langle f,f\rangle_t$, i.~e.
$\|\psi\|_t:=\sup\{|\psi(f)|:  \, \|f\|_t\le 1 \}.$
Set
\begin{equation}\label{sminus}
s_0=1-\lim_{\tau\to\infty}\sup_{t\ge 1}
\frac1{\log\tau}\log\frac {\|\psi\|_{
t}}{\|\psi\|_{\tau t}}
\end{equation}
and
\begin{equation}\label{splus}
s_1=1-\lim_{\tau\to\infty}\inf_{t\ge 1}
\frac1{\log\tau}\log\frac  {\|\psi\|_{t}}{\|\psi\|_{\tau t}}.
\end{equation}

We can specialize Theorem (\ref{interpol1}) to the this special case of
interpolating between $L^2$ and $H^1.$

\begin{Prop}\label{spinterpol} Suppose $\psi\in (H^1)^*$ and let
$Y_0=\{f\in H^1:\ \psi(f)=0\}.$  Then:\newline
(1) $(L^2,Y_0)_{s,2}=H^s$ if and only if $0\le s< s_0.$\newline
(2) $(L^2,Y_0)_{s,2}$ is a closed subspace of codimension one in $H^s$ if
and only if $s_1<s \le 1.$\end{Prop}

\begin{proof}
We  can  then  apply  Theorem  \ref{interpol1} with
$X_0=H^1 $ and $X_1=L^2$.
To estimate
$K(t,\psi)$ we  note that if $f\in H^1$ and $t\ge 1$ then
$$ \max(\|f\|_1,t\|f\|)\le \|f\|_t \le \sqrt 2\max(\|f\|_1,t\|f\|).$$
and so, for $t\ge 1,$
$$ \|\psi\|_t\le K(t^{-1},\psi)\le \sqrt 2\|\psi\|_t.$$
Hence   we   can   describe the numbers
$\sigma_0,\sigma_1$ of Theorem \ref{interpol1}  by
$$
\sigma_1=\lim_{\tau\to\infty}\sup_{t\ge 1}
\frac1{\log\tau}\log\frac {\|\psi\|_{
t}}{\|\psi\|_{\tau t}}
$$
and
$$
\sigma_0=\lim_{\tau\to\infty}\inf_{t\ge 1}
\frac1{\log\tau}\log\frac  {\|\psi\|_{
t}}{\|\psi\|_{\tau t}}.
$$
Since $\sigma_1=1-s_0$ and $\sigma_0=1-s_1$ this proves the
Proposition.\end{proof}

We next turn to the problem of estimating $\|\psi\|_t.$  The following
lemma will be useful:

\begin{Lem}\label{Fit} Suppose $F$ satisfies (\ref{type}) and
(\ref{integral}). Then for any real $t$ we have:
\begin{equation}\label{itest}
 |F(it)| \le |t|^{\frac12}e^{\pi|t|}\left(\frac1\pi
\int_{-\infty}^{\infty}\frac{|F(x)|^2}{t^2+x^2}dx\right)^{\frac12}.
\end{equation} \end{Lem}

\begin{proof} It suffices to consider $t>0.$  Then by (\ref{type2})
 $F(z)e^{i \pi z}(z+it)^{-1}$ is bounded and
analytic in the upper-half plane and so we have:
$$ F(it)= \frac{t e^{\pi t}}{\pi}\int_{-\infty}^{\infty}
\frac{F(x)}{x+it}\frac{e^{i \pi x}}{x-it}\,dx$$
Applying the Cauchy--Bunyakowski inequality  we prove the lemma.
\end{proof}

We can now give an estimate for $\|\psi\|_t$ which essentially solve
the problem of determining $s_0$ and $s_1.$

 \begin{Thm} \label{psit0} There exist a constant $C$ so
that for $t \ge 2$ we have
\begin{equation}\label{intcond} \frac1C \left(
\int_{-\infty}^{\infty}\frac{|F(x)|^2}{x^2+t^2}dx\right)^{\frac12}
\le \|\psi\|_t \le C \left(
\int_{-\infty}^{\infty}\frac{|F(x)|^2}{x^2+t^2}dx\right)^{\frac12}\end{equation}.
\end{Thm}

\begin{proof} We start with the remark that the functions
$\{(2\pi)^{-\frac12}(n^2+t^2)^{-\frac12}e_n:\ n\in\mathbb
Z\}$ together
with $(\frac12(t\sinh 2\pi t)^{-\frac12} (e_{it}+e_{-it}))$
form an orthonormal basis of $H^1$ for $\|\cdot\|_t.$  Hence
\begin{equation}\label{ortho} \|\psi\|_t^2 = \frac1{2\pi}
\sum_{n\in\mathbb
Z}\frac{|F(n)|^2}{n^2+t^2} +
4\frac{|F(it)+F(-it)|^2}{t\sinh
2\pi t}.\end{equation}
By (\ref{itest})
 the last term in
(\ref{ortho}) can be estimated by
\begin{equation}\label{ortho2}
\frac{|F(it)+F(-it)|^2}{t\sinh
2\pi t}\le C^2
\int_{-\infty}^{\infty}\frac{|F(x)|^2}{t^2+x^2}dx\end{equation}
for
$t\ge
1.$

Now if $-1\le \tau\le 1$ the map $T_{\tau}:H^1\to H^1$
defined by $T_{\tau}f=e_{\tau}f$
satisfies $\|T_{\tau}\|_t\le 2$ provided $t\le 1.$  Hence if
$\psi_{\tau}=T^*_{\tau}\psi$ we have $\frac12
\|\psi_{\tau}\|_{t}\le \|\psi\|_t\le 2\|\psi_{\tau}\|_t.$
However using (\ref{ortho}) and (\ref{ortho2}) gives:
$$\frac1{2\pi} \sum_{n\in\mathbb Z}
\frac{|F(n+\tau)|^2}{n^2+t^2}
\le \|\psi_{\tau}\|_t^2 \le
\frac1{2\pi} \sum_{n\in\mathbb Z}\frac{
|F(n+\tau)|^2}{n^2+t^2}
+ C^2 \int_{-\infty}^{\infty}
\frac{|F(x+\tau)|^2}{t^2+x^2}dx.$$
Now by integrating for $0\le\tau\le 1$ we obtain
(\ref{intcond}).
\end{proof}

\section{Application to nonharmonic Fourier series}\label{n3}
\setcounter{equation}{0}

 At this point we turn our attention to exponential Riesz bases.  Let
$\Lambda=(\lambda_n)_{n\in\mathbb Z}$ be a sequence of complex numbers.
For convenience we shall write $\sigma_n=\Re \lambda_n$ and
$\tau_n=\Im\lambda_n.$

Let us suppose that $(e_{\lambda_n})_{n\in\mathbb Z}$ is an unconditional
basis of $L^2$, or equivalently,
$((1+|\tau_n|)^{\frac12}e^{-\pi|\tau_n|}e_{\lambda_n})_{n\in\mathbb Z}$
is a Riesz basis of $L^2.$
Then this family is {\it complete interpolating set \cite{Seip95}}.
In  particular, we have
{\it sampling condition}: there exists a constant $D$ so that if
$f\in L^2$ then
\begin{equation}\label{sampling} D^{-1}\|f\|\le \left(\sum_{n\in\mathbb
Z}(1+|\tau_n|)e^{-2\pi|\tau_n|}|\hat f(\lambda_n)|^2\right)^{\frac12}\le
D\|f\| \end{equation}
(i.e., the latter family is {\it a  frame}).
We also note that it must satisfy a {\it separation condition} i.e. for
some $0<\delta<1$ we have:
\begin{equation}\label{separation}
\frac{|\lambda_m-\lambda_n|}{1+|\lambda_m-\overline{\lambda}_n|}\ge
\delta \qquad m\neq n. \end{equation}

Then we can define an entire function $F$
by
\EQ{product}{ F(z)
=\lim_{R\to\infty}\prod_{|\lambda_k|\le R}(1-\frac{z_k}{\lambda_k}).
}
We  replace  the  term  $(1-\lambda_k^{-1}z)$  by  $z$ if
$\lambda_k=0.$
We  call  $F$  the  {\it  generating function} for the unconditional
basis $(e_{\lambda_n}).$

\begin{Prop}  \cite{PW34,Levin56,I96}
The product
(\ref{product})
converges to an entire  function of exponential type  $\pi $
and   satisfies   the   integrability  conditions (\ref{integral})
and
\EQ{compl}{
\int_{-\infty }^{\infty}  \vert F(x)\vert^2 dx=\iy.
}
\end{Prop}

Let us note, that the inequality  in
(\ref{integral})  is necessary for
minimality of   the   family   and   (\ref{compl})   for  completeness  of
$(e_{\lambda_n})$.   Note that since
 $F$ satisfies (\ref{type}) and
(\ref{integral}) so that there exists $\psi\in (H^1)^*$ with
$\hat\psi=F.$  We remark that
$F$ is a Cartwright class function and then \cite{Levin56}
we have
the Blaschke condition
$$
\sum_{\lambda_n\neq 0}
\frac{|\tau_n|}{|\lambda_n|^2}<\infty
$$
Note that this implies that the families
$(e_{\lmn})_{\Im \lmn > 0}$, $(e_{\lmn})_{\Im \lmn < 0}$
are minimal in
$L^2(0,\iy)$, $L^2(-\iy,0)$ correspondingly.
Also we have
$\sum_{\lambda_n\neq
0}\frac{1}{|\lambda_n|^2}<\infty$,
(this follows from \cite{Levin} p. 127).
Thus, we have a strong Blaschke condition
\begin{equation}\label{Blaschke}
\sum_{\lambda_n\neq 0}
\frac{1+|\tau_n|}{|\lambda_n|^2}<\infty.\end{equation}

 Now by the result of Russell, Proposition
\ref{Russell},
the functions $(e_{\lambda_n})$ form an unconditional basis of a
closed  subspace  $Y_0$  of  $H^1$  of codimension one.  It is clear that
the kernel of $\psi$ coincides with $Y_0.$  Hence our above results
Proposition \ref{spinterpol} and Theorem \ref{psit0} apply
to this case.

\begin{Thm}\label{uncbasis}  Suppose $(e_{\lambda_n})_{n\in\mathbb Z}$ is
an unconditional basis of $L^2.$  Then:
\newline
(1)
$(e_{\lambda_n})_{n\in\mathbb Z}$ is an unconditional  basis
of  the  Sobolev  space   $H^s$  if  and  only   if
$0\le s< s_0.$
\newline
(2)  $(e_{\lambda_n})_{n\in\mathbb  Z}$  is an
unconditional  basis  of  a  closed  subspace  of  $H^s$  of
codimension one if and only if $s_1<s\le 1.$\newline (3)  If
$s_0\le   s\le   s_1$   then   $(e_{\lambda_n})$  is  not  an
unconditional basic sequence.  \end{Thm}

\begin{proof}
By Russell's theorem, Proposition
\ref{Russell} above,
$(e_{\lambda_n})_{n\in\mathbb Z}$ is an unconditional  basis
for a closed subspace $Y_0$ of codimension one which is  the
kernel of  the linear  functional $\psi.$  Let $v_n$  be the
weight sequence
$\ds v_n=\frac {\sinh(2\pi \tau_n)}{\tau_n}=\|e_{\lmn}\|^2_{L^2}$
and let
$h_n=(1+|\lambda_n|^2)v_n=\|e_{\lmn}\|^2_{H^1}$.
It follows from the basis property
that the map  $V:\ell_2(h)\to Y_0$ defined  by
$$
V(\alpha)=\sum_{n\in\mathbb Z} \alpha_ne_{\lambda_n}
$$
is an
isomorphism  (onto).    Clearly  $V$  is  an  isomorphism of
$\ell_2(v)$ onto $L^2(-\pi,\pi)=Y_1$. Hence by interpolation $V$ is
an isomorphism of $\ell_2(v^{1-s}h^s)$ onto  $Y_{1-s}=[Y_0,L^2]_{1-s,2}$.
In other words, setting $q_n=v_n^{1-s}h_n^s=v_n(1+|\lmn|^2)^s$, we have
$$
C^{-1} \sum |\alpha_n|^2 q_n\le \|\sum \alpha_ne_{\lmn}\|^2_{Y_{1-s}}\le
C \sum |\alpha_n|^2 q_n
$$
and the almost normalized family
$(e_{\lmn}/q_n^{1/2})_{n \in \mathbb Z}$
forms a \Rb in $Y_{1-s}$.
Thus, if $Y_{1-s}$ is a closed subspace in $H^s$, $(e_{\lambda_n})$
forms an unconditional basic sequence in $H^s$ also.

We next estimate $\|e_{\lambda_n}\|_{H^s}$ to have the inverse implication.   In fact from interpolation
between $L^2$ and $H^1$
we have
$$ \|e_{\lambda_n}\|_{H^s} \le
C\|e_{\lambda_n}\|^{1-s}\|e_{\lambda_n}\|_1^s =
C(v_n^{1-s}h_n^s)^{\frac12},
$$
where $C$ depends only on $s.$
Similarly if we define $\phi_n(f)=(f,e_{\lambda_n})$ then the norm of
$\phi_n$ in $(H^s)^*$ can be estimated by
$$
\|\phi_n\|_{(H^s)^*} \le C_1\|\phi_n\|^{1-s}\|\phi_n\|_{(H^1)^*}^s
=C_1(v_n^{1-s})^{1/2}(v_n^2/h_n)^{s/2}=C_1(v_n^{1+s}h_n^{-s})^{\frac12}.
$$
>From the other hand,
$\|e_{\lmn}\|_{H^s}\ge |\phi_n(e_{\lmn})|/\|\phi_n\|_{H^s}$,
what gives
$\|e_{\lambda_n}\|_{H^s} \ge C_1^{-1} (v_n^{1-s}h_n^s)^{1/2}$.
Therefore the norms  $\|e_{\lambda_n}\|_{H^s}$ and
$\|e_{\lambda_n}\|_{Y_{1-s}}$ are both equivalent to $\sqrt{q_n}$.
Therefore the assumption that $(e_{\lmn})$ is an unconditional basic
sequence leads to equivalence
of metrics $H^s$ and $Y_{1-s}$.
\end{proof}

{\it Remark.}
It is easy to have necessary and sufficient condition for an exponential
family $(e_{\lambda_n})_{n\in\mathbb Z}$ which is complete and minimal in
$L^2$ to be
complete and/or minimal in $H^s$ \cite{AI00}.
To do this we connect the generating function $F$ with the
critical exponent
$s_{\Lambda}$
$$
s_{\Lambda }:=\inf \lbrace s : \int_{-\infty}^{\infty}
\frac{|F(x)|^2}{1+|x|^{2s}}dx<\infty=\inf \{s: \psi\in (H^s)^* \}.
$$
Now
$(e^{i\lmn t})$  is complete in \Hs for $s<s_{\Lm}$ and is minimal for
$s>s_{\Lm}-1$. The situation for $s=s_\Lm$ or $s=s_\Lm-1$ depends
on whether $\psi$ is bounded in $H^{s_\Lm}$.
Note that $s_0\le s_{\Lambda}\le s_1$ in general.

Thus,  the family $(e_{\lmn})$ is minimal in $H^s$ for $0<s<1$,
and for any $(\a_n)\in l^2$, $\alpha\ne 0$,
$$
0<\|\sum \a_n \en/\sqrt{ q_n}\|^2_{H^s}\le
\|\sum \a_n \en /\sqrt{ q_n}\|^2_{Y_{1-s}} \le C \sum |\a_n|^2.
$$
We do not know whether $(e_{\lambda_n})_{n\in\mathbb Z}$ can be a
conditional basis  of $H^s$, for some appropriate ordering, when $s_0\le
s\le s_{\Lambda}.$
\vskip1mm

 In order to get  precise estimates  of $s_0$ and $s_1$ we will need
to
establish an alternative formula for $\|\psi\|_t$ in this special case.

Let us introduce the function $\Phi(z)$ defined by $\Phi(z)= |F(z)|
d(z,\Lambda)^{-1}$ when $z\notin \Lambda$ and
$\Phi(\lambda_n)=|F'(\lambda_n)|$ for $n\in \mathbb Z.$  The function
$\Phi$ plays an important role in the known conditions for
$(e_{\lambda_n})$ to be an unconditional basis (see \cite{Minkin} and
\cite{LS}).   We will call $\Phi$ the {\it carrier function} for
$(e_{\lambda_n}).$

The following lemma lists some useful properties:

\begin{Lem}\label{Phi} Suppose $-\infty<t<\infty.$  Then:\newline
(i) There is at most one $n\in\mathbb Z$ so that
$|it-\lambda_n|<\frac12\delta |t|$ where $\delta$ is the separation
constant in (\ref{separation}).
There is also at most one $n\in\mathbb Z$ so that
$|it-\lambda_n|<\frac14\delta|\lambda_n|.$
\newline
(ii)  $|F(it)|\le
(|\lambda_0|+|t|)\Phi(it).$ \newline
(iii) There is a constant $C$ independent of $t,n$ so that for every
$n\in\mathbb Z$ we have:
$$ |F(it)| \le
C(|\lambda_0|+|t|)\frac{|it-\lambda_n|}{(|\lambda_n|^2+t^2)^{\frac12}}\Phi(it).$$
\newline
(iv) We have  for $t\neq 0,$
$$ \Phi(it) \le |t|^{-\frac12}e^{\pi
|t|}\left(\frac1\pi\int_{-\infty}^{\infty}
\frac{|F(x)|^2}{t^2+x^2}dx\right)^{\frac12}.$$
\end{Lem}

\begin{proof} (i) Suppose $m,n$ are distinct and
$|it-\lambda_n|,|it-\lambda_m| < \frac12 \delta |t|.$  Then
$|\lambda_m-\lambda_n|<\delta t$ while
$$
|\lambda_m-\overline{\lambda_n}|
\ge |(\lambda_m-it)-\overline{(\lambda_n-it)} +2it|
\ge (2-\delta)|t|>|t|.
$$
Hence
$$
\frac{|\lambda_m-\lambda_n|}{1+|\lambda_m-\overline{\lambda}_n|} <
\delta
$$
which contradicts (\ref{separation}).

For the second part note that if
$|it-\lambda_n|<\frac14\delta|\lambda_n|$ then $|\lambda_n|<2|t|$ so that
$|it-\lambda_n|<\frac12\delta |t|.$

(ii) is immediate from the fact that $d(it,\Lambda)\le |\lambda_0|+t.$

(iii) If $|it-\lambda_n|<\frac12\delta t$ then,
in view of (i),
$|it-\lambda_n|\Phi(it)=|F(it)|$ and $t^2+|\lambda_n|^2\le 5t^2.$
Let
$|it-\lambda_n|\ge \frac12\delta t$.
Then
$$
\frac{|it-\lmn|}{(|\lmn|^2+t^2)^{1/2}}\ge
\frac{|it-\lmn|}{|\lmn|+|t|}\ge
\frac{|it-\lmn|}{|it-\lmn|+2|t|}\ge c>0.
$$
Since $|\lmn|+|t|\ge d(it,\Lm)$,
we have (iii).

(iv) Let $\lambda_n$ satisfy $|it-\lambda_n|=d(it,\Lambda).$  If $t$ and
$\tau_n$ have opposite signs or if $\tau_n=0$, then $d(it,\Lambda)\ge t$
and so that
$\Phi(it)\le t^{-1}|F(it)|.$  If they have the same sign define $G(z)=
(z-\overline{\lambda}_n)(z-\lambda_n)^{-1}F(z)$ and note that
$$ \Phi(it)= \frac{|G(it)|}{|it-\overline{\lambda}_n|} \le
t^{-1}|G(it)|.$$   Since $|G(x)|=|F(x)|$ for $x$ real, we obtain (iv)
from (\ref{itest}). \end{proof}

We next show that the Blaschke condition
(\ref{Blaschke}) can be improved for Riesz bases:

\begin{Prop}\label{Blaschke2}  If $(e_{\lambda_n})_{n\in\mathbb Z}$ is an
unconditional basis of $L^2$  then there is a
constant $C$ so that for any $0<t<\infty,$
\begin{equation}\label{Blaschke2eq} \sum_{\lambda_n\neq 0}
\frac{t(1+|\tau_n|)}{|\lambda_n|^2+t^2}\le C.\end{equation}
\end{Prop}

\begin{proof} Let us apply (\ref{sampling}) to $e_{\pm it}$.  Then
\begin{equation}\label{sampling2}
\sum_{n\in\mathbb Z}
(1+|\tau_n|)e^{-2\pi|\tau_n|}\left|\frac{\sin(\pi(\lambda_n\pm it))}
{\lambda_n \pm it}\right|^2
\le 4 D^2\|e_{\pm it}\|^2=
4 D^2 \frac{\sinh 2\pi t}{t}.
\end{equation}
Now for each $n$ there is a choice of sign so that:
$$\left|\frac{\sin(\pi(\lambda_n\pm it))}
{\lambda_n \pm it}\right| \ge \frac{|\sinh
(\pi(|\tau_n|+t))|}{|\lambda_n|+t}$$ and hence
$$ \sum_{n\in\mathbb Z}
(1+|\tau_n|)e^{-2\pi|\tau_n|}\frac{\sinh(\pi(t+|\tau_n|)^2}{|\lambda_n|^2+t^2}
\le  4D^2  \frac{\sinh 2\pi t}{t}.$$
This yields (\ref{Blaschke2eq}) for $t\ge 1$ and this extends to $t\ge
0$ in view of (\ref{Blaschke}) and the fact that
$\sum_{n\neq 0}|\lambda_n|^{-2}<\infty.$\end{proof}

We will also need a perturbation lemma:

\begin{Lem}\label{perturbation} Let $(e_{\lambda_n})$ and
$(e_{\mu_n})_{n\in\mathbb Z}$ be two unconditional bases of $L^2.$
  Suppose further that there is a constant
$C$ so that
$$ \sum_{n\in\mathbb Z}
\frac{t|\mu_n-\lambda_n|}{|\mu_n||\lambda_n|+t^2}
\le C
\qquad 1<t<\infty.$$
Suppose $\Phi$ and $\Psi$ are the carrier functions for
$(e_{\lambda_n})$ and
$(e_{\mu_n}).$
Then there exist constants
$B,T>0$ so that if
$t\ge
T$
$$ \frac1B\frac{\Psi(it)}{\Phi(it)} \le \prod_{\substack{0<|\lambda_n|\le
t\\|\mu_n|\neq 0}}\frac{|\lambda_n|}{|\mu_n|} \le B \frac
{\Psi(it)}{\Phi(it)}.$$
\end{Lem}

\begin{proof}
We observe that for each
$n$ we have (taking $t=\max(1,|\mu_n|^{\frac12}|\lambda_n|^{\frac12})$)
$$ |\lambda_n-\mu_n| \le C\max(1,
2|\mu_n|^{\frac12}|\lambda_n|^{\frac12}).$$    Hence
$$
|\lmn|\le |\mu_n| +2C|\lmn|^{1/2}|\mu_n|^{1/2} +C\le |\mu_n|+
\frac12|\lambda_n|+ 2C^2|\mu_n| +C,$$
Along with a similar estimate
for $|\mu_n|$ and setting $C_1=2+4C^2>1$ we get:
\begin{equation}\label{Csub1}
 |\lambda_n| \le C_1(|\mu_n|+1), \quad
 |\mu_n|\le C_1(|\lambda_n|+1).\end{equation}

 Now let
$c=\frac14\min(\delta,\delta')$ where $\delta,\delta'$ are the separation
constants of $(\lambda_n)_{n\in\mathbb Z}$ and $(\mu_n)_{n\in\mathbb Z}$
respectively.

 We next make the remark that there is a constant $M$ so that if
$|w|,|z|\le 2C_1+1$ and $|1-w|,|1-z|\ge c$ then
\begin{equation}\label{loglips} |\log|1-w|-\log|1-z||\le
M|w-z|.\end{equation}

Let us fix $T=|\mu_0|+|\lambda_0|+2C_1.$  Suppose that $t\ge T,$
and   let
$p=p(t),q=q(t)\in\mathbb Z$ be chosen so that
$|it-\lambda_p|=\min\{|it-\lambda_n|: \ n\in\mathbb Z\}$ and
$|it-\mu_q|=\min\{|it-\mu_n|:\ n\in\mathbb Z\}.$  It may happen that
$p=q.$  Note that we have an automatic estimate,
\EQ{2t}{ |it-\lambda_p|\le |t|+|\lambda_0|\le 2t,\quad |it-\mu_q|\le
|t|+|\mu_0|\le 2t.
}
Then if $n\neq p,q$ and $|\lambda_n|>t$ we have
$|\mu_n|>\frac12C_1^{-1}|\lambda_n|$ and so
$|it-\mu_n|\le |t|+|\mu_n|\le (2C_1+1)|\mu_n|.$
 By Lemma \ref{Phi} (i) we have
$|it-\lambda_n|\ge c|\lambda_n|$ and $|it-\mu_n|\ge c|\mu_n|.$  Hence we
have by (\ref{loglips})
\begin{align*}
|\log|it-\mu_n|-\log|it-\lambda_n|-\log|\mu_n|+\log|\lambda_n| |&\le
M\frac{t|\lambda_n-\mu_n|}{|\lambda_n||\mu_n|}\\ &\le (2C_1+1)M
\frac{t|\lambda_n-\mu_n|}{|\lambda_n||\mu_n|+t^2}.\end{align*}
Next suppose $n\neq p,q$ and $|\lambda_n|\le t.$  Then $|\mu_n|\le
C_1(t+1) \le 2C_1t.$  We also have $|it-\lambda_n|,|it-\mu_n|\ge c|t|$
and so by (\ref{loglips})
\begin{align*}
|\log|it-\mu_n|-\log|it-\lambda_n| |&\le
M\frac{|\lambda_n-\mu_n|}{t}\\ &\le (2C_1+1)M
\frac{t|\lambda_n-\mu_n|}{|\lambda_n||\mu_n|+t^2}.\end{align*}
Combining and summing over all $n\neq p,q$ we have
$$\log\frac{\Psi(it)}{\Phi(it)}=
\delta(t)\log\left|\frac{it-\mu_p}{it-\lambda_q}\right|
+\sum_{0<|\lambda_n|\le t}\log|\lambda_n|
-\sum_{\substack{0<|\lambda_n|\le t\\ |\mu_n|\neq 0}}\log|\mu_n|
+\gamma(t)$$
where $|\gamma(t)|\le C(2C_1+1)M$ and $\delta(t)=1$ if $p\neq q$ and $0$
if $p=q.$

To conclude we need only consider the case $p\neq q.$  In this case
$|it-\mu_p|,|it-\lambda_q|\ge ct$.  We also have $|\lambda_p|,|\mu_q| \le
3t$ by
(\ref{2t}) and so by (\ref{Csub1}) $|\lambda_q|,|\mu_p|\le C_1(3t+1)\le
4C_1t.$ Hence
$|it-\mu_p|,|it-\lambda_q|\le 5C_1t.$  This concludes the
proof.\end{proof}

\begin{Lem}\label{plusminus} Suppose $(e_{\lambda_n})$ is an
unconditional basis of $L^2.$  Then there exist constants $B,T$ so that
if $t\ge T$ then
$$\frac1B \Phi(it)\le \Phi(-it)\le B\Phi(it).$$\end{Lem}

\begin{proof}  This follows from Lemma \ref{perturbation}  taking
$\mu_n=\overline{\lambda}_n$ in view of Lemma \ref{Blaschke2}.\end{proof}

The next Theorem is the key step in the proof of our main result:

\begin{Thm}\label{psitnormthm} Suppose $(e_{\lambda_n})_{n\in\mathbb Z}$
is an unconditional
basis of $L^2$.  Then there is a constant $C$ and $T>0$
so that if $t\ge T$ then
\EQ{psitnorm}{
C^{-1}t^{\frac12} e^{-\pi t}\Phi(it)\le
\|\psi\|_t
\le
Ct^{\frac12}e^{-\pi t}\Phi(it). }
\end{Thm}

\begin{proof}The left-hand inequality in (\ref{psitnorm}) is an
immediate consequence of Lemma \ref{Phi} (iv) and
(\ref{intcond}). We turn to the right-hand inequality.

We first use Lemma \ref{Phi}(ii), (iii) and Lemma \ref{plusminus}.
There are constants
$C,T>1$ so that if
$|t|\ge T$ we have $\Phi(-it)\le C\Phi(it)$, $|F(it)|\le
Ct\Phi(it)$ and for every
$n,$
\begin{equation}\label{FPhiest} |F( it)| \le Ct\frac{|\lambda_n-
it|}{(|\lambda_n|^2+t^2)^{\frac12}} \Phi( it).\end{equation}

 Choose  $g\in H^1$ so
that $\psi(f)=\langle f,g\rangle_t$ for $f\in H^1.$ Let  $h$ be the
orthogonal projection with respect to $\langle\cdot\rangle_t$ of $g$
onto the subspace $H^1_0$ of all $f$ so that $f(-\pi)=f(\pi)=0$ and let
$k=g-h.$ Then $\|\psi\|_t^2=\|k\|_t^2+\|h\|_t^2.$

The orthogonal complement of $H^1_0$ (with respect to
$\langle\cdot\rangle_t$) is a 2-dimensional space with orthonormal basis
$\{e_{\pm it}/\|e_{\pm it}\|_t\}.$  Hence
$$k=\|e_{it}\|_t^{-2}(\overline{F(it)}e_{it}
+\overline{F(-it)}e_{-it})$$ and
$$ \|k\|_t^2 = \|e_{it}\|_t^{-2}(|F(it)|^2+|F(-it)|^2).$$  Since
$\|e_{it}\|_t^2 =2t\sinh 2\pi t$, we deduce
$$ \|k\|_t \le C_1t^{\frac12} \Phi(it)e^{-\pi t}\qquad t\ge
T$$ and a suitable constant $C_1.$
 It remains therefore only to estimate
$\|h\|_t.$

We first argue that
\begin{align*} \langle e_z,k\rangle_t &=(2t\sinh 2\pi t)^{-1}(
F(it)\langle e_z,e_{it}\rangle_t
+ F(-it)
\langle e_z,e_{-it}\rangle_t)\\
&= \frac{i}{\sinh 2\pi t}
(F(-it)\sin\pi(z-it)-F(it)\sin\pi (z+it)).\end{align*}
Since $\psi(e_{\lambda_n})=F(\lambda_n)=0$ for $n\in\mathbb Z$ we deduce
that
$$ \langle e_{\lambda_n},h\rangle_t = \frac{i}{\sinh2\pi t}
(F(it)\sin\pi(\lambda_n+it)-F(-it)\sin\pi (\lambda_n-it)).$$

 Now if we use (\ref{FPhiest} we  get
an estimate valid for $t\ge T:$
 \begin{equation*} |\langle e_{\lambda_n},h\rangle_t|
\le
C\Phi(it)\frac{t|\lambda_n+it||\lambda_n-it|}
{(|\lambda_n|^2+t^2)^{\frac12}\sinh 2\pi t}
\left(\left|\frac{\sin
(\pi(\lambda_n-it))}{\lambda_n-it}\right|+\left|\frac
{\sin(\pi(\lambda_n+it))}{\lambda_n+it}\right|\right).\end{equation*}

Since $h\in H^1_0$ we then have $$\langle e_{\lambda_n},h \rangle_t
=(\lambda_n^2+t^2) (e_{\lambda_n},h) $$
and we can then rewrite the above estimate as
$$
|(e_{\lambda_n},h)| \le C\frac{\Phi(it)}{\sinh 2\pi t}
\frac{1}{(|\lmn|^2 +t^2)^{1/2}}
\left(\left|\frac{\sin
(\pi(\lambda_n-it))}{\lambda_n-it}\right|+\left|\frac
{\sin(\pi(\lambda_n+it))}{\lambda_n+it}\right|\right).$$

Now
$$  (e_{\lambda_n},th+h') =
(t-i\lambda_n)(e_{\lambda_n},h).
$$

 We next use the sampling inequality
(\ref{sampling}):
$$
 \|h\|_t^2 =
\|th+h'\|_{L^2}^2 \le D^2 \sum_{n\in\mathbb
Z}(1+|\tau_n|)e^{-2\pi|\tau_n|}
|t- i\lmn |
|(e_{\lambda_n},h)|^2.
$$
However we can combine with (\ref{sampling2}) to
deduce that
$$ \|th+h'\|_{L^2} \le 4C^2D^2 t^{\frac12}\Phi(it)(\sinh 2\pi
t)^{-\frac12}$$ for
$t\ge T$ which gives the conclusion.\end{proof}

We now consider the case when $(\lambda_n)$ is a small perturbation of
the sequence $\mu_n=n.$  For convenience we shall assume that
$\lambda_n=0$ can only occur when $n=0.$

\begin{Thm} \label{perturbation2}
Suppose  $(e_{\lambda_n})_{n\in\mathbb Z}$ is an unconditional basis of
$L^2$ and for some constant $C$ and all $t\ge 1$
\begin{equation}\label{perturb} \sum_{n\neq 0}
\frac{t|\lambda_n-n|}{n^2+t^2}<C.\end{equation}
Then
\begin{equation*}\label{splus1}  s_1=\frac12 + \lim_{\tau\to\infty}
\sup_{t\ge
1}\frac{1}{\log\tau}\sum_{t< |n|\le\tau
t}\log\frac{|n|}{|\lambda_n|}\end{equation*} and
\begin{equation*}\label{sminus1}
s_0=\frac12 + \lim_{\tau\to\infty}
\inf_{t\ge
1}\frac{1}{\log\tau}\sum_{t< |n|\le\tau
t}\log\frac{|n|}{|\lambda_n|}.\end{equation*}  \end{Thm}

\begin{proof} In this case we compare the carrier function $\Phi$ for the
basis $(e_{\lambda_n})$ with the carrier function $\Psi$ for the
basis
$(e_n).$  Clearly $\Psi(it)=|\sin\pi it|/\pi t.$  We can next use Lemma
\ref{perturbation} to estimate $\Phi(it)$  and then the theorem follows
directly from
   Theorem \ref{psitnormthm}
together with (\ref{splus}) and (\ref{sminus}).  \end{proof}

Let us specialize to some important cases.   Let
$\delta_n=\Re\lambda_n-n=\sigma_n-n.$

\begin{Thm}\label{perturbation3}
Suppose  $(e_{\lambda_n})_{n\in\mathbb Z}$ is an unconditional basis of
$L^2$ such that $\sup|\delta_n|<\infty$ and $\sum_{n\neq
0}\tau_n^2n^{-2}<\infty.$  Then
\begin{equation}\label{splus2}  s_1=\frac12 - \lim_{\tau\to\infty}
\inf_{t\ge
1}\frac{1}{\log\tau}\sum_{t< |n|\le\tau
t}\frac{\delta_n}{n}\end{equation} and
\begin{equation}\label{sminus2}
s_0=\frac12 - \lim_{\tau\to\infty}
\sup_{t\ge
1}\frac{1}{\log\tau}\sum_{t< |n|\le\tau
t}\frac{\delta_n}{n}.\end{equation}  \end{Thm}

{\it Remark.}  In particular (\ref{splus2}) and (\ref{sminus2}) hold if
 $|\lambda_n-n|$ is bounded.

\begin{proof}
Combining
Proposition \ref{Blaschke2} and the boundedness of $(\delta_n)$ gives us
(\ref{perturb}).  Note that if $n\neq 0,$
$$ \log \frac{|n|}{|\lambda_n|}= -\log (1
+\frac{|\lambda_n|-|n|}{|n|}).$$
Now
\begin{align*} \frac{|\lambda_n|-|n|}{|n|} &= \left(1+
\frac{2\delta_n}{n}+\frac{\delta_n^2+\tau_n^2}{n^2}\right)^{\frac12} \\
&= \frac{\delta_n}{n}+\alpha_n\end{align*} where
$$ |\alpha_n| \le C\frac{1+\tau_n^2}{n^2}$$ for a suitable constant $C.$
By (\ref{Blaschke})  and the assumption of the theorem,
this implies that $\sum_{n\neq 0}|\alpha_n|<\infty$ and
yields the Theorem.\end{proof}

Before discussing examples we observe one more  property of $s_0$ and
$s_1$ in this case, which uses recent results of \cite{LS} and the theory
of $A_2-$weights.

\begin{Thm}\label{strictineq}  If $(e_{\lambda_n})_{n\in\mathbb Z}$ is an
unconditional basis of $L^2$ then $s_0>0$ and $s_1<1.$\end{Thm}

\begin{proof}  We will use the connections between Riesz basis property and
sampling/interpolation in the spaces of entire functions of
exponential type. These connections in the case of $L^2$ and the
Paley--Wiener space may be found in \cite{Seip95}.

Let $L^2_{\pi,s}$ consisting of all
entire functions of exponential type at most $\pi$ and satisfying
$$ \int_{-\infty}^{\infty}
\frac{|f(\xi)|^2}{(1+|\xi|)^{2s}}d\xi<\infty.$$
(Note that the Fourier transform of $L^2_{\pi,s}$ is
the set of all distributions from $H^{-s}(\mathbb R)$ supported on
$[-\pi,\pi].$)  Now the formal adjoint of the map from $\ell_2(\mathbb
Z)$ to $ H^s$ defined by
$(\alpha_n)\to
\sum_{n\in\mathbb Z}\alpha_n
(1+|\tau_n|)^{-\frac12}(1+|\lambda_n|)^{-s}e_{\lambda_n}$
is  the map from $L^2_{\pi,s} $ to $\ell_2(\mathbb Z)$ given  by
$f\to
\l(f(\lambda_n)(1+|\tau_n|)^{1/2}
(1+|\lambda_n|)^{s}e^{-\pi\tau_n}\r)_{n\in\mathbb Z}.$  Hence
$(e_{\lambda_n})_{n\in\mathbb Z}$ is an unconditional basic sequence
(resp.
unconditional basis) if and only if $(\lambda_n)_{n\in\mathbb Z}$ is an
interpolating sequence (resp. complete interpolating sequence) in $
L^2_{\pi,s}.$

Note that if $(\lambda_n)$ is
interpolating for $L^2_{\pi,s-1}$ then it is interpolating for
$L^2_{\pi,s}$ by the simple device of considering functions of the form
$f(z)= (z-\mu)g(z)$ where $\mu\notin\Lambda=\{\lambda_n\}_{n\in\mathbb
Z}$ and
$g\in L^2_{\pi,s-1}.$

It follows that our result can be proved by showing that
$(\lambda_n)_{n\in\mathbb Z}$ is a complete interpolating sequence for
$L^2_{\pi,s}$ for all $|s|<\epsilon$ for some $\epsilon>0.$  To do this
we use the results of \cite{LS} that  this is equivalent to requiring
that $(1+|\xi|)^{2s}\Phi(\xi)^2$
 is an
$A_2$-weight for $|s|<\epsilon.$ Now $\Phi^2$ is an $A_2$-weight
(\cite{LS} or \cite{Minkin})
and so there exists $\eta>0$ so that $\Phi^{2(1+\eta)}$ is an
$A_2-$weight
(cf.
\cite{Garnett} p. 262, Corollary 6.10).  Hence the Hilbert transform is
bounded on both $L_2(\mathbb R, \Phi^{2(1+\eta)})$ on $L^2(\mathbb R,
(1+|\xi|)^{2\th})$ for $0<\th<\frac12.$ It then follows by
complex interpolation
that
$\Phi(\xi)^2(1+|\xi|)^{2s}$ is an $A_2-$weight when
$|s|<\eta(1+\eta)^{-1}.$\end{proof}

Note that these results now imply Theorem \ref{main}.

{\bf Examples.}
We recall  the classical  theorem of  Kadets, see, e.g.,
\cite{KNP} or \cite{Levin},   that   if   $(\lambda_n)$   are   real
then
$\sup_n|\delta_n|<\frac14$  is  a  sufficient  condition for
$(e_{\lambda_n})_{n\in\mathbb  Z}$  to  be  a  Riesz  basis.
First,  we  consider  the  case  of  regular  behavior.  For
example, we  can set  $\dl_n=-\frac{1}{2}q\operatorname{sign}
n$, see \cite{Avd74}.   Then we obtain $s_1=s_0=\frac12+q.$ More
generally if for some $y>0$ we have $C^{-1}(1+|x|)^{2q}\le |F(x+iy)|\le
C(1+|x|)^{2q}$, we obtain $s_1=s_0 = \frac12+q$
(if we use the integral estimates of $\| \psi\|_t$, i.e. Theorem
\ref{psit0}).

One    can    easily    make    sequences   $(\dl_n)$   with
$\sup|\dl_n|<\frac14$ to exhibit any required behavior.   In
fact if we  put $$ b_n  = \frac{1}{\log 2}  \sum_{2^n<|k|\le
2^{n+1}}\frac{\dl_k}{k}$$         then         $$        s_0
=\frac12-\lim_{N\to\infty}\frac1N                 \inf_{n\ge
1}\sum_{k=n+1}^{n+N}b_k$$     and     $$     s_1=    \frac12
-\lim_{N\to\infty}\frac1N                         \sup_{n\ge
1}\sum_{k=n+1}^{n+N}b_k.$$

To be  more specific if $-\frac12<p\le q<\frac12$,  set
\begin{eqnarray*}
&\dl_n=\frac12q\operatorname{sign} n\
\
\
&\text{for} \ 2^{(2^{2k})}< |n| \le 2^{(2^{2k+1})}\\
&\dl_n=\frac12p\operatorname{sign} n\
\
\
&\text{for} \ 2^{(2^{2k-1})}< |n| \le 2^{(2^{2k})}  \end{eqnarray*}
Then      for      $      2^{2k}<m\le      2^{2k+1}$      $$
b_m=\frac{q}{\log2}\sum_{2^m+1}^{2^{m+1}}\frac   1k=q+o(1)
$$ and for $ 2^{2k-1}<m\le 2^{2k}$ $$ b_m=p+o(1).  $$ Thus
$$  s_0=\frac12-  q,\  s_1=\frac12-  p  $$  (note that an example
of  irregular behavior is given  in \cite {Avd74})


\begin{thebibliography}{10}



\bibitem{Avd74}
S.A. Avdonin,
\newblock On {R}iesz bases of exponentials in $L^2$,
\newblock {\em Vestnik Leningrad Univ., Ser. Mat., Mekh., Astron.},(1974)
(13):5--12,
\newblock (Russian); English transl. in Vestnik Leningrad Univ. Math., v. 7
  (1979), 203-211.

\bibitem{AI95}
S.~Avdonin and S.A. Ivanov,
\newblock {\em Families of Exponentials. The Method of Moments in
  Controllability Problems for Distributed Parameter Systems},
\newblock Cambridge University Press, Cambridge, 1995.

\bibitem{AI00}
S.~Avdonin and S.~Ivanov.
\newblock {L}evin--{G}olovin theorem for the {S}obolev spaces.
\newblock {\em Matemat. Zametki}, 68(2):{163--172}, 2000.
\newblock (Russian); English transl. in Math. Notes.




\bibitem{AIR}
S.~A. Avdonin, S.~A. Ivanov, and D.~L. Russell,
\newblock Exponential bases in {S}obolev spaces in control and observation
  problems for the wave equation,
\newblock {\em Proceedings of the Royal Society of Edinburgh},
\newblock to appear.

\bibitem{BS} C. Bennett and R. Sharpley, {\em Interpolation of
operators,} Academic Press, Orlando, 1988.

\bibitem{Garnett} J.B. Garnett, {\em Bounded analytic functions},
Academic Press, Orlando, 1981.

\bibitem{I96}
S.A. Ivanov,
\newblock Nonharmonic {F}ourier series in the {S}obolev spaces of positive
  fractional orders,
\newblock {\em New Zealand Journal of Mathematics}, 25 (1996) 36--46,

\bibitem{IP78}
S.A. Ivanov and B.S. Pavlov,
\newblock Carleson series of resonances in the {R}egge problem,
\newblock {\em Izv. Akad. Nauk SSSR, Ser. Mat.}, 42(1) (1978) 26--55,
\newblock (Russian); English transl. in Math. USSR Izvestija., v. 12, no. 1
  (1978), 21--51.

\bibitem{Janson}
S.~Janson,
\newblock Interpolation of subcouples and quotient couples,
\newblock {\em Ark. Mat.}, 31 (1993) 306--338.

\bibitem{KS} S. Kaijser and P. Sunehag, Interpolation of subspaces of
codimension one, in preparation.

\bibitem{KNP}
S.~V. Khrushchev, N.~K. Nikol'skii, and B.~S. Pavlov,
\newblock {\em Unconditional bases of exponentials and reproducing kernels},
  volume 864 of {\em Lecture Notes in Math.``Complex Analysis and Spectral
  Theory''}.
\newblock Springer--Verlag, Berlin/Heidelberg, 1981.
\newblock pp. 214--335.

\bibitem{KMP} N. Krugliak, L. Maligranda and L.E. Persson,\newblock The
failure
of the Hardy inequailty and interpolation of intersections,\newblock Ark.
Mat.
37
(1999) 323--344.




\bibitem{LLT}
I.~Lasiecka, J.-L. Lions, and R.~Triggiani,
\newblock Nonhomogeneous boundary value problems for second order hyperbolic
  operators,
\newblock {\em J. Math. Pures et Appl.}, 65(2) (1986) 149--192.


\bibitem{Levin56}
B.Ya Levin,
\newblock {\em Distribution of zeros of entire functions},
\newblock GITTL, Moscow, 1956,
\newblock En\-glish transl., Amer.Math.Soc., Prov\-i\-dence, RI,1964.

\bibitem{Levin} B. Ya. Levin, {\em Lectures on entire functions,}
Translations of Mathematical Monographs, Vol. 150, Amer. Math. Soc.,
Providence 1996.

\bibitem{LM}
J.-L. Lions and E.~Magenes,
\newblock {\em Probl\`emes aux limites nonhomog\'enes et applications}, volume
  I,II.
\newblock Dunod, Paris, 1968.

\bibitem{L1} J. L\"ofstr\"om, \newblock Real interpolation with
constraints, \newblock J. Approx. Theory
82
(1995), 30--53.

\bibitem{L2} J. L\"ofstr\"om, \newblock Interpolation of subspaces,
\newblock Preprint, University of Goteborg, 1997 \newblock
http://www.math.chalmers.se/\~{}jorgen.


\bibitem{LS}
Yu. Lyubarskii and K.~Seip,
\newblock Weighted Paley-Wiener spaces,
\newblock  submitted


\bibitem{Minkin}
A.M. Minkin,
\newblock Reflection of exponents, and unconditional bases of exponentials,
\newblock {\em Algebra i Anal.}, 3(5) (1992) 110--135.
\newblock (Russian); English transl. in St. Petersburg Math. J., v. 3, no. 5
  (1992), 1043--1068.


\bibitem{NS86}
K.~Narukawa and T.~Suzuki,
\newblock Nonharmonic {F}ourier series and its applications,
\newblock {\em Appl. Math. Optim.}, 14 (1986) 249--264.

\bibitem{Nik}
N.~K. Nikol'ski\u\i,
\newblock {\em A Treatise on the Shift Operator},
\newblock Springer, Berlin, 1986.

\bibitem{PW34}
R.E.A.C. Paley and N.~Wiener,
\newblock {\em Fourier Transforms in the Complex Domain}, volume XIX,
\newblock AMS Coll. Publ., New York, 1934.

\bibitem{Russell82}
D.L. Russell,
\newblock On exponential bases for the {S}obolev spaces over an
interval,
\newblock {\em J. Math. Anal. Appl.}, 87(2) (1982) 528--550.

\bibitem{Seip95}
K.~Seip,
\newblock On the connection between exponential bases and certain related
  sequences in {$L^2(-\pi,\pi)$},
\newblock {\em J. of Functional Analysis}, 130(1) (1995) 131--160.


\bibitem{Wallsten}
R.~Wallsten,
\newblock Remarks on interpolation of subspaces,
\newblock {\em Lecture Notes in Math.}, 1302 (1988) 410--419.

\end{thebibliography}
\end{document}